\documentclass[reqno,twoside]{amsart}
\usepackage{amsfonts}
\usepackage{amsmath,amssymb,colordvi}
\usepackage{multicol}
\usepackage{color,bbm}
\usepackage{epic}
\usepackage{graphicx}
\usepackage{xcolor}
\usepackage{esvect}
\usepackage{cite}
\usepackage{amsmath,amssymb,amsfonts,latexsym,amsthm}
\usepackage{algorithmic}
\usepackage{graphicx}
\usepackage{algorithm,algorithmic}
\usepackage{hyperref}
\usepackage{enumerate}

\usepackage{a4wide,amsmath,amssymb,latexsym,amsthm}
\usepackage{graphicx}

\newtheorem{theorem}{Theorem}[section]
\newtheorem{proposition}[theorem]{Proposition}
\newtheorem{lemma}[theorem]{Lemma}
 
\newtheorem{definition}[theorem]{Definition}
\newtheorem{remark}[theorem]{Remark}
\theoremstyle{definition}

\usepackage{float}

\newtheorem{example}{Example}
\numberwithin{equation}{section}

\graphicspath{{./Imgs/}}

\newcommand{\R}{\mathbb R}
 
\newcommand{\N}{\mathbb N}

\newcommand{\be}{\begin{equation}}
\newcommand{\ee}{\end{equation}}
\newcommand{\ba}{\begin{eqnarray}}
\newcommand{\ea}{\end{eqnarray}}
\newcommand{\beq}{\begin{equation}}
\newcommand{\eeq}{\end{equation}}

\usepackage{color}

\usepackage[textsize=small]{todonotes}
\numberwithin{equation}{section}

\def\N{{\mathbb{N}}}

\usepackage{color}
\usepackage{stmaryrd}

\subjclass[2010]{}

\begin{document}

\title[Tracking control]{Tracking controllability for finite-dimensional linear systems}

\author{Sebasti\'an Zamorano}
\address[S.  Zamorano]{University of Deusto, Faculty of Engineering, Av. Universidades 24, 48007, Bilbao, Basque Country, Spain.}
 \email{sebastian.zamorano@deusto.es}

 \author{Enrique Zuazua}
\address[E. Zuazua]{Friedrich-Alexander-Universität Erlangen-Nürnberg, Department of Data Science, Chair for Dynamics, Control and Numerics (Alexander von Humboldt Professorship), Cauerstr. 11, 91058 Erlangen, Germany.
\newline \indent   Chair of Computational Mathematics, Fundación Deusto,	Avenida de las Universidades, 24, 48007 Bilbao, Basque Country, Spain.
\newline \indent  Universidad Autónoma de Madrid, Departamento de Matemáticas, Ciudad Universitaria de Cantoblanco, 28049 Madrid, Spain.}
\email{enrique.zuazua@fau.de}

\begin{abstract}
This paper develops a functional-analytic characterization of output tracking controllability for finite-dimensional linear systems. By formulating tracking as the surjectivity of the control-to-output map on suitable trajectory spaces, we show that exact tracking is equivalent to a trajectory-space observability inequality associated with the dual input-output structure. This characterization enables a Hilbert Uniqueness Method (HUM) type variational construction of minimum-norm tracking controls and makes explicit the intrinsic regularity requirements on reference trajectories induced by the system dynamics and the output operator. The same framework also yields a natural notion of approximate tracking when exact tracking fails. We provide explicit formulas in the scalar case and report illustrative numerical examples for ODEs and semi-discretized PDEs, demonstrating the method for both smooth and nonsmooth targets.
\end{abstract}

\maketitle


\section{Introduction, problem formulation and main results}

The problem of forcing the output of a dynamical system to follow a
prescribed reference trajectory, usually referred to as output tracking,
reference tracking, or trajectory planning, is a classical topic in
systems and control theory. It arises naturally in applications such as
robotics, aerospace engineering, process control and autonomous systems,
where the objective is not only to steer the state, but to reproduce a
desired time evolution of a measured or prescribed output.

In this paper, we consider the finite-dimensional, time-invariant linear
input-output systems of the form
\begin{equation}\label{eq:0}
\begin{cases}
x'(t)=Ax(t)+Bu(t) & t\in(0,T),\\
x(0)=x_0,
\end{cases}
\end{equation}
with output
\begin{equation}
\label{eq:0.1}
y(t)=Cx(t).
\end{equation}
Here $T>0$ is fixed, \(x:[0,T]\to \mathbb{R}^n\) denotes the state,
\(u:[0,T]\to \mathbb{R}^m\) is the control input, \(x_0\in\mathbb{R}^n\),
and $A\in\mathbb{R}^{n\times n}$, $B\in\mathbb{R}^{n\times m}$, $C\in\mathbb{R}^{p\times n}$, with $1\leq m\leq n$ and $1\leq p\leq n$.

The exact output tracking problem consists in prescribing a reference
trajectory \(f\) and looking for a control $u$ such that the
corresponding solution of \eqref{eq:0} satisfies $Cx(t)=f(t)$, $t\in [0,T]$. More precisely, given a reference trajectory $f$ in a target space $\mathcal Y_C$, to be specified according to the input-output structure of the triple $(A,B,C)$, and satisfying the corresponding compatibility conditions at $t=0$, we seek $u\in L^2(0,T;\mathbb{R}^m)$ such that the solution $x$ of \eqref{eq:0} fulfills
\begin{equation}\label{eq:0.2}
Cx(t)=f(t),\qquad t\in [0,T].
\end{equation}
Accordingly, throughout the paper the tracking property is understood as
a property of the input-output triple \((A,B,C)\), or equivalently of the
system \eqref{eq:0}-\eqref{eq:0.1}, and not of the state equation \eqref{eq:0} alone.

Tracking problems for finite-dimensional linear systems have a long
history and have been studied from several complementary viewpoints.
The early formulation of reproducibility for multivariable systems goes
back to Brockett and Mesarovi\'c~\cite{BrockettMesarovic1965}, and was
later developed through the theory of invertibility of linear systems,
right-invertibility and inverse systems~\cite{Silverman1969,SainMassey69,Ros70}.
Closely related notions of path controllability and target path
controllability were also studied in the linear input-output and
time-varying settings~\cite{AlbrechtGrasseWax1986,Wohltmann1985}.
Classical algebraic and geometric approaches rely on notions such as
invertibility, invariant zeros, relative degree and polynomial matrix
methods; see, for instance,~\cite{Ros70,SilvermanPayne1971,IsidoriBook}.
These tools provide structural solvability criteria and are closely
related to dynamic inversion and exact tracking by differentiation of the
desired output. In particular, the relative degree identifies how many
time derivatives of the output are needed before the control appears
explicitly, and therefore gives precise regularity requirements on
admissible reference trajectories. In the scalar-output case with relative degree $r$, the natural choice
is $\mathcal Y_C=H^r_0(0,T),$ since $r$ is the number of derivatives of the output needed before the
control appears explicitly. For instance, the choice $\mathcal Y_C=H^1_0(0,T;\mathbb{R}^p)$ corresponds to outputs of
relative degree one, or to a simplified functional setting in which only first-order trajectory regularity is retained.

Related notions such as functional controllability, functional observability and output $\varepsilon$-controllability, together with
earlier works on path controllability \cite{AlbrechtGrasseWax1986}, \cite{Wohltmann1985} and target path controllability, also focus on reachable output trajectories and are therefore close in spirit to the present formulation  \cite{DavronLissy2025,garcia2013alternative,GermaniMonaco83,darouach2025functional}. Differential flatness \cite{Fliess95}  provides a constructive parametrization of states and controls through suitable flat outputs, but it is not, in general, a theory of tracking a fixed prescribed output $Cx$.

Against this background, the present work should be viewed as a complementary variational approach to output tracking. Its main contribution is not the trajectory-space viewpoint itself, which already appears in classical input-output controllability theories, but the use of the Hilbert Uniqueness Method (HUM) to analyze output tracking for the triple $(A,B,C)$. In this context, HUM means that the tracking control is constructed through a dual variational problem posed on the adjoint equation, rather than by directly inverting the input-output dynamics. More precisely, exact tracking is recast as the surjectivity of the control-to-output map between Hilbert trajectory spaces and is characterized by a dual observability inequality for the adjoint equation with source $C^\mathsf{T} g$ and observation $B^\mathsf{T}\phi$. This leads to a trajectory-space Gramian acting on the admissible output space $\mathcal Y_C$, and to a variational construction of the minimum-$L^2$-norm tracking control. Compared with explicit inversion methods, this approach separates the regularity and compatibility requirements encoded in $\mathcal Y_C$ from the construction of the control itself, and selects the minimum-$L^2$-norm tracking control. Moreover, when exact tracking is not available, the same duality framework yields a natural penalized formulation for approximate tracking. This is useful already in finite dimensions and is particularly suited to semi-discrete PDE settings, where explicit algebraic inversion may become ill-conditioned or unavailable.

Therefore, the main contributions of the paper are the following:
\begin{enumerate}
\item We formulate exact output tracking controllability for the input-output system \eqref{eq:0}--\eqref{eq:0.1} as the surjectivity of the control-to-output operator
$$
\Lambda_C:L^2(0,T;\mathbb R^m)\longrightarrow \mathcal Y_C,
$$
where the admissible trajectory space $\mathcal Y_C$ encodes the regularity and compatibility constraints imposed by the triple
$(A,B,C)$. This surjectivity is characterized through a dual observability inequality and an associated tracking Gramian acting on
trajectory spaces.

\item We develop a HUM-type construction of minimum-norm tracking controls through a dual variational problem posed on the adjoint equation endowed with the observation $B^\mathsf{T}\varphi$. The resulting control is expressed in terms of the minimizer of a quadratic functional and the corresponding adjoint solution.

\item We introduce approximate output tracking controllability in $L^2(0,T;\mathbb R^p)$ and prove its equivalence with a unique
continuation property for the adjoint equation endowed with the observation $B^\mathsf{T}\varphi$.

\item In the scalar-input case, we derive explicit formulas in controller canonical coordinates. These formulas recover the classical relative-degree regularity requirements and show how compatibility conditions among several output components arise from the system dynamics.

\item We present illustrative numerical experiments for ODE systems and semi-discrete PDE models, showing how the variational formulation behaves for smooth and nonsmooth targets.
\end{enumerate}

%
%
%
%

The paper is organized as follows. Section~\ref{sec:2} introduces the
control-to-output operator and proves the equivalence between exact
output tracking, surjectivity, dual observability and invertibility of a
tracking Gramian. Section~\ref{sec:3} develops the HUM formulation and the
minimum-norm tracking control. Section~\ref{sec:4} studies approximate output
tracking and its equivalence with a unique continuation property. In
Section~\ref{sec:5} we analyze the scalar-input case in controller canonical form,
highlighting relative degree and compatibility conditions. Section~\ref{sec:6}
contains numerical experiments, and Section~\ref{sec:8} concludes with comments
and open problems.


\subsection{Notation}
Let $\mathcal Y_C$ be a Hilbert space of admissible output trajectories. We assume that $\mathcal Y_C$ is continuously embedded in $L^2(0,T;\mathbb R^p)$, and we denote its dual by $\mathcal Y_C'$. The duality pairing between $\mathcal Y_C'$ and $\mathcal Y_C$ will be denoted by $\langle \cdot,\cdot\rangle_{\mathcal Y_C',\mathcal Y_C}$. For $r\in\N$, we use the notation $H_0^r(0,T):=H^r_{0,\mathrm{left}}(0,T):=\{\phi\in H^r(0,T): \ \phi^{(k)}(0)=0, \ k=0,\ldots,r-1\}$. No condition at $t=T$ is imposed. The $L^2$ inner product is $\langle \xi,\phi\rangle_{L^2}=\int_0^T \xi^\mathsf{T}(t)\phi(t)dt$. Matrix transposition is denoted by $\mathsf{T}$. 
\section{Tracking controllability: first results}\label{sec:2}
In this section, we recast the exact output tracking problem in an operator-theoretic framework that allows one to make precise the duality between tracking controllability and observability properties of the adjoint system. This formulation provides the natural setting for the variational constructions developed later via the HUM. For simplicity and without loss of generality, we restrict our analysis to zero initial conditions $x_0=0$ and targets $f$ satisfying  $f(0)=0$. This simplifies the presentation while preserving the essential structural features of the tracking problem. The general case with nonzero initial data can be recovered by a standard translation argument $\tilde f(t)=f(t)-Ce^{At}x_0$, $\tilde y(t)=Cx(t)-Ce^{At}x_0$.

\subsection{The input–output operator}
For $x_0=0$, the output tracking condition $Cx(t)=f(t)$ becomes
\begin{equation*}
f(t)=\int_0^t Ce^{A(t-s)}Bu(s)\,ds.
\end{equation*}
Motivated by this identity, we define the control-to-output operator
$\Lambda_C:L^2(0,T;\mathbb{R}^m)\to \mathcal Y_C$ by
\begin{equation*}
(\Lambda_Cu)(t)=\int_0^t Ce^{A(t-s)}Bu(s)\,ds.
\end{equation*}
Rather than viewing tracking as a pointwise condition on trajectories, we interpret it as a mapping property between spaces of controls and output trajectories, in the spirit of functional/output controllability and path-controllability approaches
\cite{GermaniMonaco83,AlbrechtGrasseWax1986,Wohltmann1985}.

Throughout the sequel, $\mathcal{Y}_C$ is fixed a priori according to the regularity and compatibility restrictions induced by the input-output triple $(A,B,C)$, and is chosen so that $\Lambda_C: L^2(0,T; \mathbb{R}^m)\to \mathcal{Y}_C$ is well-defined and continuous. It is not defined as $\operatorname{Ran}(\Lambda_C)$. Hence exact output tracking controllability is a surjectivity property relative to the prescribed admissible output space $\mathcal Y_C$.
\begin{definition}\label{def:2}
The input-output system \eqref{eq:0}-\eqref{eq:0.1} is said to be exactly output tracking controllable if the operator $\Lambda_C$ is surjective. 
\end{definition}

Characterizing surjectivity of $\Lambda_C$ naturally leads us to its adjoint, anticipating observability inequalities and revealing which output components can be effectively reconstructed from the control action. To this end, we examine $\Lambda_C^*$.
\begin{lemma}
The adjoint operator $\Lambda_C^*: \mathcal{Y}'_C \to L^2(0,T;\mathbb{R}^m)$ is given by $\Lambda_C^*\psi=B^{\mathsf{T}}\varphi_\psi,$
where $\varphi_\psi$ is the transposition solution of $-\varphi'(t)=A^{\mathsf{T}}\varphi(t)+C^{\mathsf{T}}\psi(t)$, $t\in(0,T),$
 and $\varphi(T)=0$. In particular, if $\psi\in L^2(0,T;\mathbb R^p)$, then
\begin{equation}\label{lambda_C^*}
(\Lambda_C^* \psi)(t) = \int_t^T B^\mathsf{T} e^{A^\mathsf{T}(\tau-t)}C^\mathsf{T} \psi(\tau) d\tau.
\end{equation}
\end{lemma}


\begin{proof}
For $L^2$-sources, \eqref{lambda_C^*} follows from Fubini's theorem. For $\psi\in\mathcal Y_C'$, $C^\mathsf{T}\psi$ is defined by transposition through $\langle C^\mathsf{T}\psi,v\rangle=\langle\psi,Cv\rangle_{\mathcal Y_C',\mathcal Y_C}$. Testing the state equation against the corresponding transposition solution $\varphi_\psi$,  we obtain $\Lambda_C^*\psi=B^T\varphi_\psi$. See, for instance, \cite[Ch.~4]{Rudinbook}
\end{proof}

\subsection{The tracking Gramian operator}
We define the tracking Gramian operator by $G_C:=\Lambda_C\Lambda_C^*:\mathcal Y_C'\rightarrow \mathcal Y_C.$ Equivalently, $G_C$ is characterized by the duality identity $\langle h,G_C g\rangle_{\mathcal Y_C',\mathcal Y_C} = \langle \Lambda_C^* g,\Lambda_C^* h\rangle_{L^2(0,T;\mathbb R^m)},$ $g,h\in\mathcal Y_C'.$ Thus $G_C$ is the operator associated with the symmetric bilinear form $ a(g,h):= \langle \Lambda_C^* g,\Lambda_C^* h\rangle_{L^2(0,T;\mathbb R^m)}.$ Unlike the classical finite-dimensional controllability Gramian, $G_C$ acts on trajectory spaces and encodes the dual observability structure of the output tracking problem. Throughout the main results we assume, without loss of generality, that $C$ has full row rank. Otherwise, one first replaces the output space by $\operatorname{Ran}(C)$ and projects the target accordingly.
\begin{lemma}
\label{lem:operator-equivalence}
The following operator-theoretic properties are equivalent:
\begin{enumerate}[$\roman{enumi})$]
    \item\label{it:i} $\Lambda_C: L^2(0,T;\mathbb{R}^m) \to \mathcal{Y}_C$ is surjective.
    \item $\Lambda_C^*: \mathcal{Y}'_C \to L^2(0,T;\mathbb{R}^m)$ is bounded below, i.e., there exists $\gamma > 0$ such that
    \begin{equation}\label{ecuacion5new}
    \|g\|_{\mathcal{Y}'_C}\leq \gamma\|\Lambda_C^* g\|_{L^2(0,T;\mathbb{R}^m)}, \quad \forall g\in \mathcal{Y}'_C.
    \end{equation}
    \item The tracking Gramian $G_C = \Lambda_C \Lambda_C^*: \mathcal{Y}'_C \to \mathcal{Y}_C$ is an isomorphism.
\end{enumerate}
Moreover, if any of these holds, the $L^2$-minimum-norm control achieving $\Lambda_C u = f$ is given explicitly by $u_{\min} = \Lambda_C^* G_C^{-1} f.$
\end{lemma}

\begin{proof}
The equivalence between $i)$ and $ii)$ follows from the closed range theorem. Indeed, $\Lambda_C$ is surjective if and only if $\Lambda_C^*$ is bounded below. 

Assume ii). Define $ a(g,h):=\langle \Lambda_C^*g,\Lambda_C^*h\rangle_{L^2(0,T;\mathbb R^m)}$,  $g,h\in\mathcal Y_C'.$ This bilinear form is continuous and, by ii), coercive on $\mathcal Y_C'$. Hence, by the Lax--Milgram theorem, for every $f\in\mathcal Y_C$ there exists a unique $g\in\mathcal Y_C'$ such that  $a(g,h)=\langle h,f\rangle_{\mathcal Y_C',\mathcal Y_C}$,  $h\in\mathcal Y_C'.$ Since  $a(g,h)= \langle h,\Lambda_C\Lambda_C^*g\rangle_{\mathcal Y_C',\mathcal Y_C} = \langle h,G_Cg\rangle_{\mathcal Y_C',\mathcal Y_C},$ we obtain $G_Cg=f$. Thus $G_C:\mathcal Y_C'\to\mathcal Y_C$ is an isomorphism. Conversely, if $G_C$ is an isomorphism and $\Lambda_C^*g=0$, then $G_Cg=0$, hence $g=0$. Moreover, $\|\Lambda_C^*g\|_{L^2}^2 = \langle g,G_Cg\rangle_{\mathcal Y_C',\mathcal Y_C}.$ Since $G_Cg=\Lambda_C\Lambda_C^*g$, the boundedness of $\Lambda_C$ gives $\|G_Cg\|_{\mathcal Y_C}\leq \|\Lambda_C\|\|\Lambda_C^*g\|_{L^2}$. If $G_C^{-1}$ is bounded, then $\|g\|_{\mathcal Y_C'}\leq \|G_C^{-1}\|\|G_Cg\|_{\mathcal Y_C}\leq \|G_C^{-1}\|\|\Lambda_C\|\|\Lambda_C^*g\|_{L^2}.$
\end{proof}

%

The following theorem provides several equivalent characterizations of output tracking controllability, linking operator surjectivity, dual observability, and Gramian invertibility.
\begin{theorem}\label{th:1}
The following assertions are all equivalent:
\begin{enumerate}[$\Roman{enumi})$]
\item The input-output system \eqref{eq:0}-\eqref{eq:0.1} is exactly output tracking controllable in the sense of Definition \ref{def:2}.
\item The operator $\Lambda_C: L^2(0,T;\mathbb{R}^m)\to \mathcal{Y}_C$ is surjective.
\item There exists a constant $\gamma>0$ such that the observability inequality holds
\begin{equation*}
\|\psi\|_{\mathcal{Y}'_C}\leq \gamma\|\Lambda_C^* \psi\|_{L^2(0,T;\mathbb{R}^m)},   \quad \forall \psi\in \mathcal{Y}'_C
\end{equation*}
\item The tracking Gramian operator $G_C: \mathcal{Y}'_C\to \mathcal{Y}_C$ is invertible with continuous inverse.
\end{enumerate}
The equivalences follow directly from Lemma \ref{lem:operator-equivalence} together with the definitions of tracking controllability and the observability inequality.
\end{theorem}

\begin{proof}
By definition,  $I)$  and  $II)$ are equivalent.  Given that $i)$ and $ii)$ in Lemma \ref{lem:operator-equivalence} are equivalent, and that $ii)$ ensures that $\Lambda_C^*$ is bounded-below,  $II)$ is equivalent to $III)$.  The remaining equivalences follow from Lemma \ref{lem:operator-equivalence}.\end{proof}
\section{Hilbert Uniqueness Method for tracking control}\label{sec:3}

The operator-theoretic characterization developed in Section \ref{sec:2} shows that exact output tracking is equivalent to a quantitative observability inequality for the adjoint equation endowed with the observation $B^T\phi$, thus reflecting the full triple $(A,B,C)$. This section investigates such observability properties and clarifies the structural mechanisms that make exact tracking possible through the use of the HUM \cite{lions1988controlabilite}, a variational technique originally developed for PDE control.
\subsection{Output tracking observability}
Let $g\in\mathcal Y_C'$. The adjoint equation 
\begin{equation} \label{eq:adjoint} 
\begin{cases} 
-\varphi'(t)=A^T\varphi(t)+C^Tg(t), & t\in(0,T),\\ 
\varphi(T)=0, 
\end{cases} 
\end{equation} 
is understood in the sense of transposition. More precisely, the source term $C^Tg$ is defined by duality through $\langle C^Tg,v\rangle := \langle g,Cv\rangle_{\mathcal Y_C',\mathcal Y_C},$ for every admissible state test function $v$ such that $Cv\in\mathcal Y_C$ and $\varphi$ satisfies
\begin{multline}\label{eq:variational}
\int_0^T \langle v'(t),\varphi(t)\rangle_{\mathbb R^n}\,dt = \int_0^T \langle Av(t),\varphi(t)\rangle_{\mathbb R^n}\,dt \\+ \langle g,Cv\rangle_{\mathcal Y_C',\mathcal Y_C}.
\end{multline}

The connection to output tracking controllability becomes clear upon noting that, for $\psi=g$, the output of the adjoint system is precisely  $B^\mathsf{T}\varphi$, which coincides with $\Lambda_C^* g$ as given in \eqref{lambda_C^*}. 

In contrast with classical state observability, the tracking problem requires estimating the source term in the adjoint system from partial observations of the adjoint trajectory. This motivates the following notion of tracking observability with respect to the observation $B^\mathsf{T}\varphi$.

\begin{definition}\label{def:obs}
We say that the input-output triple $(A,B,C)$ satisfies the output tracking
observability inequality on $[0,T]$ if there exists $\gamma>0$ such
that, for every $g\in \mathcal{Y}'_C$, the solution $\varphi$
of \eqref{eq:adjoint} satisfies
\begin{equation}\label{eq:obs}
\|g\|_{\mathcal{Y}'_C}^2\leq \gamma\|B^\mathsf{T} \varphi\|_{L^2(0,T;\mathbb{R}^m)}^2.
\end{equation}
Inequality \eqref{eq:obs} will be referred to as the output tracking observability inequality.
\end{definition}

\begin{remark}
Even though $g$ acts on the adjoint system through $C^T g$, the inequality in Definition \ref{def:obs} is stated in terms of $g$ itself; the resulting structural restriction to $\operatorname{Ran}(C)$ is addressed in Section \ref{sec:3:2}.
\end{remark}

The following theorem shows that the output tracking observability inequality for the triple $(A,B,C)$ implies exact output tracking
controllability and provides a variational construction of the minimum-norm control.
\begin{theorem}\label{th:2}
Assume that the input-output triple $(A,B,C)$ satisfies the output tracking
observability inequality. Then, for every target $f\in \mathcal{Y}_C$, there exists  a unique control $u_{\min}\in L^2(0,T;\mathbb{R}^m)$ of minimal $L^2$-norm such that $\Lambda_C u_{\min}=f$. Moreover, $u_{\min}$ is given by
\begin{equation}\label{eq:control:min}
u_{\min}(t)=B^\mathsf{T}\widehat{\varphi}(t),
\end{equation}
where $\widehat{\varphi}$ is the solution of \eqref{eq:adjoint} corresponding to the unique minimizer $\widehat{g}$ of the quadratic functional
\begin{equation}\label{eq:functional}
J(g):=\frac{1}{2}\int_0^T |B^\mathsf{T}\varphi(t)|^2 dt -\langle g,f\rangle_{\mathcal{Y}'_C,\mathcal{Y}_C}, \ \forall g\in \mathcal{Y}'_C.
\end{equation}  
\end{theorem}

\begin{proof}
We proceed in three steps.

{\bf Step 1:} Let $x$ be the solution of \eqref{eq:0} with $x(0)=0$ and control $u$.  Using $x$ as a test function in the weak formulation \eqref{eq:variational} with $v=x$, and  $\langle C^\mathsf{T}g,x\rangle_{\mathcal{Y}'_C,\mathcal{Y}_C}= \langle g, Cx\rangle_{\mathcal{Y}'_C,\mathcal{Y}_C}$, we obtain, after integration by parts, that the output tracking condition $Cx=f$ is equivalent to
\begin{equation}\label{eq:critical-points}
\langle g,f\rangle_{\mathcal{Y}'_C,\mathcal{Y}_C}=\int_0^T \langle B^\mathsf{T} \varphi(t),u(t)\rangle_{\mathbb{R}^m}dt,\quad \forall g\in \mathcal{Y}'_C,
\end{equation}
where $\varphi\in L^2(0,T;\mathbb{R}^n)$ solves \eqref{eq:adjoint} with source $g$.\\

{\bf Step 2:} Consider the quadratic functional $J$ defined in \eqref{eq:functional}. Its G\^ateaux derivative at $\widehat{g}$ in the direction $g$ is
\begin{align*}
d J(\widehat{g};g)=\int_{0}^T\langle B^\mathsf{T}\widehat{\varphi}(t), B^\mathsf{T}\varphi(t)\rangle_{\mathbb{R}^m} dt - \langle g,f\rangle_{\mathcal{Y}'_C,\mathcal{Y}_C},
\end{align*}
where $\widehat{\varphi}$ and $\varphi$ correspond to $\widehat{g}$ and $g$, respectively. Setting $dJ(\widehat{g};g)=0$ for all $g$ yields precisely condition \eqref{eq:critical-points} with $u=B^\mathsf{T}\widehat{\varphi}$.  Thus, any critical point $\widehat{g}$ of $J$ provides a control $u=B^\mathsf{T}\widehat{\varphi}$ achieving exact output tracking.\\

{\bf Step 3:}
The functional $J$ is convex and continuous on $\mathcal{Y}'_C$. The output tracking observability inequality \eqref{eq:obs} implies coercivity. Indeed, we observe that
\begin{align*}
J(g)\geq \frac{1}{2\gamma}\|g\|_{\mathcal{Y}'_C}^2 
 - \|g\|_{\mathcal{Y}'_C}\|f\|_{\mathcal{Y}_C},
\end{align*}
which tends to infinity as $\|g\|_{\mathcal{Y}'_C} \to \infty$. Hence, by the direct method of the calculus of variations, $J$ attains a unique minimizer $\widehat{g}\in \mathcal{Y}'_C$. The corresponding control $u_{\min}=B^\mathsf{T}\widehat{\varphi}$ satisfies \eqref{eq:critical-points} and is the unique control of minimal  $L^2$-norm achieving $\Lambda_C u=f$. \end{proof}

\subsection{Tracking observability inequality: preliminary results}\label{sec:3:2}

In the present subsection, we restrict ourselves to the model case $\mathcal Y_C=H^1_0(0,T;\mathbb R^p)$, $ \mathcal Y_C'=H^{-1}(0,T;\mathbb R^p).$ This choice corresponds, for instance, to outputs of relative degree one, and allows us to state explicit preliminary estimates. The purpose of the following result is not to give a full characterization of \eqref{eq:obs} for arbitrary target spaces $\mathcal Y_C$, but rather to illustrate how source terms in the adjoint equation may be estimated from the observation $B^T\varphi$ under suitable structural assumptions.

Throughout this subsection we assume that the pair $(A,B)$ is
controllable, or equivalently that the pair $(A^T,B^T)$ is observable,
in the Kalman sense $\operatorname{rank}\,[\,B,AB,A^2B,\ldots,A^{n-1}B\,]=n$, see \cite{kalman1960}.
Under this assumption, the homogeneous adjoint dynamics satisfies the
classical finite-dimensional observability inequality with observation $B^T\varphi$.

\begin{proposition}\label{th:obs}
   Assume that $(A,B)$ is controllable.  Let $g(t)=\sum_{k=0}^d g_k t^k$, where $g_k\in \R^p$, for every $k\in\{0,1,\ldots,d\}$. Let $\varphi\in L^2(0,T;\R^n)$ be the solution by transposition of \eqref{eq:adjoint}. Then, there exists a constant $\eta>0$, depending only on $A,B,C,T$ and $d$, such that
    \begin{equation}\label{eq:16}
        \|g\|_{H^{-1}(0,T;\R^p)}\leq \eta\|B^\mathsf{T}\varphi\|_{H^{d+1}(0,T;\R^m)}.
    \end{equation}
   If $C$ does not have full row rank, the estimate holds with $g$ replaced by $\mathbb P_{\operatorname{Ran}(C)}g$ on the left-hand side, where $\mathbb P_{\operatorname{Ran}(C)}$ is the orthogonal projection onto $\operatorname{Ran}(C)$.
\end{proposition}

\begin{proof} (Sketch)
We begin with the simplest situation where $g$ is constant in time $(k=0)$.
Let $g(t)\equiv g\in\R^p$, for every $t\in(0,T)$, and $\varphi\in L^2(0,T;\mathbb{R}^n)$ be the unique  weak solution \eqref{eq:adjoint}.
Differentiating \eqref{eq:adjoint} with respect to time, the derivative $\psi=\varphi'$ satisfies
\begin{equation*}
    \begin{cases}
        -\psi'(t)=A^\mathsf{T}\psi(t) & t\in (0,T),\\
        \psi(T)=-C^\mathsf{T} g.
    \end{cases}
\end{equation*}
Thus, by classical observability estimates for ODEs, there exists a constant $c>0$ such that:
\begin{equation*}
   \|C^\mathsf{T} g\|_{\R^n}^2 \leq c\int_0^T|B^\mathsf{T} \psi(t)|^2 dt =c\|B^\mathsf{T}\varphi\|_{H^1(0,T;\R^m)}^2.
\end{equation*}

For $d\geq 1$, we proceed by induction. The key idea is to differentiate the adjoint equation repeatedly, each time reducing the polynomial degree of the source term until a constant source is reached. At each step, observability inequality applied to the differentiated system yields an estimate for a coefficient $g_k$ in terms of higher derivatives of $B^\mathsf{T}\varphi$. Summation over $k$ and the continuity of the embedding of polynomials into $H^{-1}$ yield \eqref{eq:16}.
\end{proof}

\begin{remark}
The estimate in Proposition~\ref{th:obs} controls only the component of the source term $g$ that effectively acts on the adjoint dynamics. Source terms taking values in $\ker(C^T)$ lie in the kernel of the adjoint operator $\Lambda_C^*$ and are therefore completely invisible to the observation $B^\mathsf{T}\varphi$. Such components correspond to output directions that are intrinsically untrackable. Consequently, output tracking observability can only be expected on the subspace of sources in $\operatorname{Ran}(C)$. This limitation is structural and independent of the regularity of the source term. This characterization identifies $\operatorname{Ran}(C)$ as the maximal class of output directions for which exact output tracking observability can be expected, a viewpoint that will be exploited in Sections~\ref{sec:4} and~\ref{sec:5}.
\end{remark}

\section{Approximate output tracking controllability}\label{sec:4}
Exact output tracking depends on the choice of the target space $\mathcal Y_C$,
which encodes the regularity imposed by the relative degree. Nevertheless, the variational HUM framework remains meaningful and naturally leads to a relaxed notion of tracking, in which the objective is to approximate a desired output in the $L^2$-sense. We now formalize this notion of approximate output tracking controllability.

\begin{definition}
The input-output system \eqref{eq:0}-\eqref{eq:0.1} is said to be approximately output tracking controllable if,
for any output target $f\in L^2(0,T;\mathbb{R}^p)$ and every $\varepsilon>0$,
there exists a control $u\in L^2(0,T;\mathbb{R}^m)$ such that the corresponding
output satisfies
\begin{equation}\label{approximate}
\|Cx(t)-f(t)\|_{L^2(0,T;\mathbb{R}^p)} \le \varepsilon.
\end{equation}
\end{definition}

\begin{definition}
The input-output triple $(A,B,C)$ is said to satisfy the unique continuation property if, whenever $g\in L^2(0,T;\mathbb{R}^p)$ and $\varphi$ is the corresponding solution of \eqref{eq:adjoint}, the following condition holds
\begin{equation}\label{eq:ucp}
B^\mathsf{T} \varphi(t) =0 \quad \text{for a.e. }t\in [0,T], \text{ implies } g\equiv 0.
\end{equation}
\end{definition}

The following result is the qualitative counterpart of the exact tracking characterization obtained above. While exact output tracking is equivalent to a quantitative observability inequality for the adjoint equation, approximate output tracking only requires the corresponding qualitative property: unique continuation from the observation $B^\mathsf{T}\varphi$. Thus, the passage from exact to approximate tracking mirrors the classical duality between observability estimates and unique continuation principles.
\begin{theorem}\label{th:approx}
The input-output system \eqref{eq:0}-\eqref{eq:0.1} is approximately output tracking controllable if and only if the triple $(A,B,C)$ satisfies the unique continuation property \eqref{eq:ucp}.
\end{theorem}

\begin{proof} The proof is based on a convex variational formulation. Although the functional is not Fr\'echet differentiable at $g=0$, convexity guarantees a well-defined subdifferential optimality condition.

We proceed with the proof in several steps.

{\bf Step 1:} Assume first that the unique continuation property \eqref{eq:ucp} holds. For a given $f\in L^2(0,T;\mathbb{R}^p)$ and $\varepsilon>0$, consider the functional $J_{\varepsilon}:L^2(0,T;\mathbb{R}^p)\to \mathbb{R} $ defined by  
\begin{equation*}
J_\epsilon(g):=\frac{1}{2}\int_0^T |B^\mathsf{T}\varphi(t)|^2 dt -\langle g,f\rangle_{L^2(0,T;\mathbb{R}^p)}+\varepsilon\|g\|_{L^2(0,T;\mathbb{R}^p)},
\end{equation*}
where $\varphi$ solves \eqref{eq:adjoint} with source $g$.  
Immediately, $J_\epsilon$ is a continuous and convex functional. We show it is coercive. 

Indeed, let $\{g_k\}_{k\geq 1}\subset L^2(0,T;\mathbb{R}^p)\setminus\{0\}$ be a sequence of source terms for the adjoint system \eqref{eq:adjoint} such that $\|g_k\|_{L^2(0,T;\mathbb{R}^p)}\to \infty$ as $k\to\infty$. Let $\varphi_k$ be the corresponding solution of \eqref{eq:adjoint} and $\widetilde{\varphi}_k\in C([0,T];\mathbb{R}^n)$ be the corresponding solution with source term $\widetilde{g}_k=g_k/\|g_k\|_{L^2(0,T;\mathbb{R}^p)}$. Then,
\begin{equation}\label{eq:aprox-funct}
\frac{J_{\varepsilon}(g_k)}{\|g_k\|_{L^2(0,T;\mathbb{R}^p)}}=\frac{1}{2}\|g_k\|_{L^2(0,T;\mathbb{R}^p)}\int_0^T |B^\mathsf{T} \widetilde{\varphi}_k(t)|^2 dt + \varepsilon - \langle \widetilde{g}_k, f\rangle_{L^2(0,T;\mathbb{R}^p)}.
\end{equation}

We observe that if $$\liminf_{k \to \infty}\int_0^T |B^\mathsf{T} \widetilde{\varphi}_k(t)|^2 dt>0,$$ the coercivity follows directly from \eqref{eq:aprox-funct}. Otherwise, let us consider the most delicate situation where we assume that $\int_0^T |B^\mathsf{T} \widetilde{\varphi}_k(t)|^2 dt \to 0$. Since $\{\widetilde{g}_{k}\}_{k}$ is bounded in $L^2(0,T;\mathbb{R}^p)$, the corresponding solution $\{\widetilde{\varphi}_{k}\}_{k}$ of the adjoint system \eqref{eq:adjoint} are uniformly bounded in $H^1(0,T;\mathbb{R}^n)$. Therefore, extracting subsequences  there exists $z\in H^1(0,T;\mathbb{R}^n)$ such that $\widetilde{\varphi}_{k_j}\rightharpoonup z$ weakly in $H^1(0,T;\mathbb{R}^n)$. By the Rellich--Kondrachov compactness theorem, $\widetilde{\varphi}_{k_j}\rightarrow z$ strongly in $L^2(0,T;\mathbb{R}^n)$. Moreover, since $\widetilde{g}_{k_j}\rightharpoonup  h$ weakly in $L^2(0,T;\mathbb{R}^p)$, passing to the limit in the adjoint equation yields that $z$ is the unique weak solution of
\begin{equation*}
\begin{cases}
-z'(t)=A^\mathsf{T}z(t)+ C^\mathsf{T} h(t) & t\in(0,T),\\
z(T)=0.
\end{cases}
\end{equation*}
We also have
\begin{equation*}
\int_0^T |B^\mathsf{T}z(t)|^2dt\leq \liminf_{j\to\infty}\int_0^T |B^\mathsf{T}\widetilde{\varphi}_{k_j}(t)|^2 dt =0.
\end{equation*}
That is, $B^\mathsf{T}z =0$ for a.e. $t\in(0,T)$. By \eqref{eq:ucp} we deduce that $h\equiv 0$. Accordingly, $\widetilde{g}_{k_j}\rightharpoonup  0$ and then $\langle\widetilde{g}_{k_j},f\rangle_{L^2(0,T;\mathbb{R}^p)}\to 0$, as $j\to\infty$. From \eqref{eq:aprox-funct}, we obtain
\begin{equation*}
\liminf_{j\to\infty}\frac{J_{\varepsilon}(g_{k_j})}{\|g_{k_j}\|_{L^2(0,T;\mathbb{R}^p)}}\geq \varepsilon
\end{equation*}
which implies $J_{\varepsilon}(g_k)\to\infty$ as  $\|g_k\|_{L^2(0,T;\mathbb{R}^p)}$ tends to infinity, and thus $J_{\varepsilon}$ is coercive and attains a minimizer $\overline{g}\in L^2(0,T;\mathbb{R}^p)$.

{\bf Step 2:} Let $\varphi_g$ be the solution of \eqref{eq:adjoint} with source $g$. For any $w\in L^2(0,T;\mathbb{R}^p)$ with corresponding solution $\varphi_w$, optimality of $g$ yields, after dividing by $|\delta|$ and letting $\delta\to 0^{\pm}$,
\begin{equation}\label{eq:aprox-final}
\left|\int_0^T\Big(\langle B^\mathsf{T}\varphi_g(t), B^\mathsf{T}\varphi_w(t)\rangle_{\mathbb{R}^m}- \langle w(t),f(t)\rangle_{\mathbb{R}^p}\Big)dt\right|\leq \varepsilon\|w\|_{L^2(0,T;\mathbb{R}^p)}.
\end{equation}

{\bf Step 3:} Set  $u=B^\mathsf{T}\varphi_g$. Using $u$ in \eqref{eq:0} with $x$, replacing the identity
into \eqref{eq:aprox-final}, we finally obtain
\begin{equation*}
\left|\int_0^T\langle Cx(t)-f(t),w(t)\rangle_{\mathbb{R}^p}dt\right|\leq \varepsilon\|w\|_{L^2(0,T;\mathbb{R}^p)},
\end{equation*}
that is, \eqref{approximate} holds.

{\bf Step 4:} Conversely, assume approximate output tracking controllability holds, and let $g\in L^2(0,T;\mathbb{R}^p)$ satisfy $B^\mathsf{T}\varphi=0$ a.e.,  $\varphi$ solving \eqref{eq:adjoint}. For any control $u\in L^2(0,T;\mathbb{R}^m)$, testing \eqref{eq:0} against $\varphi$ yields
\begin{equation*}
\int_0^T \langle Cx(t),g(t)\rangle_{\mathbb{R}^p}dt=\int_0^T \langle B^\mathsf{T}\varphi(t),u(t)\rangle_{\mathbb{R}^m} dt=0.
\end{equation*}
However, given that the system fulfills the property of approximate output tracking controllability, it follows that the subspace of the projected trajectories of the form $Cx(t)$ is dense in $L^2(0,T;\mathbb{R}^p)$. This implies that $g \equiv 0$. \end{proof}

\begin{remark}
Theorem \ref{th:approx} shows that approximate output tracking is governed by a qualitative unique continuation property of the adjoint system, rather than by a quantitative observability inequality. While this guarantees density of the reachable outputs in $L^2$, it does not provide information on control costs or exact reachability, which require the stronger observability estimates developed in Sections~\ref{sec:2}--\ref{sec:3}.
\end{remark}

\section{Scalar-input case: controller form and compatibility conditions}\label{sec:5}

In this section we specialize the previous abstract results to the
scalar-input case. The purpose is twofold. First, we recall explicitly
the controller canonical form and the corresponding coordinate change.
Second, we use this form to identify the regularity requirements and
compatibility conditions that arise in exact output tracking.

Let $x'(t)=Ax(t)+bu(t)$, $x(0)=0,$ where $A\in\mathbb{R}^{n\times n}$, $b\in\mathbb{R}^n$, and
$u:[0,T]\to\mathbb{R}$. We assume throughout this section that the pair
$(A,b)$ is controllable, namely $\operatorname{rank}\,[\,b,Ab,\ldots,A^{n-1}b\,]=n$. Then there exists an invertible matrix $P\in\mathbb{R}^{n\times n}$ such that, in the new coordinates $x=Pz,$ the system takes the controller canonical form $z'(t)=\widetilde A z(t)+\widetilde b\,u(t),$ where $\widetilde A=P^{-1}AP$, $\widetilde b=P^{-1}b.$

With the convention used here, we take
\begin{equation*}
\widetilde A=
\begin{pmatrix}
0 & 1 & 0 & \cdots & 0\\
0 & 0 & 1 & \cdots & 0\\
\vdots & \vdots & \vdots & \ddots & \vdots\\
0 & 0 & 0 & \cdots & 1\\
-a_0 & -a_1 & -a_2 & \cdots & -a_{n-1}
\end{pmatrix},
\qquad
\widetilde b=
\begin{pmatrix}
0\\
0\\
\vdots\\
0\\
1
\end{pmatrix}.
\end{equation*}
The output matrix is transformed accordingly $y(t)=Cx(t)=CPz(t)=\widetilde C z(t)$, $\widetilde C:=CP.$ Thus the output tracking problem for the original triple $(A,b,C)$ is equivalent to the output tracking problem for the transformed triple
$(\widetilde A,\widetilde b,\widetilde C)$.

\subsection{Single-output tracking and relative degree}

We first consider the case of a scalar output, $y(t)=Cx(t),$  $C\in\mathbb{R}^{1\times n}.$ Writing $\widetilde C=CP=(\eta_1,\ldots,\eta_n)$, we define $k_*:=\max\{j\in\{1,\ldots,n\}:\eta_j\neq 0\}.$ The integer $r:=n-k_*+1$ is the relative degree of the output $y=Cx$. Equivalently, $r$ is the smallest positive integer such that $CA^{r-1}b\neq 0$. In controller coordinates, this number determines how many derivatives of the prescribed output are needed before the control appears explicitly.

\begin{theorem}\label{th:v1}
Assume that $(A,b)$ is controllable and let $C\in\mathbb{R}^{1\times n}$. Let $r$ be the relative degree of the output $y=Cx$. Then, for every target $f\in H^r_0(0,T)$, there exists a control $u\in L^2(0,T)$ such that the corresponding solution satisfies $Cx(t)=f(t)$, $t\in[0,T].$ Moreover, in controller coordinates the control can be obtained by differentiating the relation $\widetilde C z=f$ exactly $r$ times.
\end{theorem}

\begin{proof}
Since the system is in controller canonical form, the first $n-1$ components satisfy $z_1'=z_2$, $z_2'=z_3$, $\ldots$,  $z_{n-1}'=z_n$, while the last equation contains the control: $z_n'=-a_0z_1-a_1z_2-\cdots-a_{n-1}z_n+u.$
The output is $y(t)=\widetilde C z(t)=\sum_{j=1}^n \eta_j z_j(t).$
By definition of $k_*$, the last nonzero coefficient of $\widetilde C$ is $\eta_{k_*}$. Hence the control does not appear in
$y,y',\ldots,y^{(r-1)}$, while it appears linearly in $y^{(r)}$ with
coefficient $CA^{r-1}b\neq 0$. Therefore, after differentiating the
tracking condition $y(t)=f(t)$ $r$ times, one obtains an identity of the form $f^{(r)}(t)
= CA^r x(t)+CA^{r-1}b\,u(t).$ Since $CA^{r-1}b\neq0$, the control is determined by
\begin{equation*}
u(t)
=
\frac{1}{CA^{r-1}b}
\left(
f^{(r)}(t)-CA^r x(t)
\right).
\end{equation*}
Thus, if $f\in H^r_0(0,T)$, the right-hand side belongs to $L^2(0,T)$, and the resulting control produces $Cx=f$. The compatibility conditions at \(t=0\) are encoded in the space $H^r_0$ under the zero initial condition.
\end{proof}

\begin{remark}
This result recovers, in the present tracking framework, the classical regularity requirement associated with relative degree. If the output has relative degree $r$, exact tracking by an $L^2$-control requires the reference trajectory to have $r$ derivatives in $L^2$. The HUM formulation of Sections~\ref{sec:2}--\ref{sec:3} gives an equivalent dual interpretation of this fact in terms of an observability inequality for the adjoint tracking system.
\end{remark}

\subsection{Several outputs and compatibility conditions}

We now consider the case of several prescribed outputs,
\begin{equation*}
y(t)=Cx(t)\in\mathbb{R}^p,
\quad
C=
\begin{pmatrix}
C_1\\
\vdots\\
C_p
\end{pmatrix},
\quad \widetilde C=CP=
\begin{pmatrix}
\widetilde C_1\\
\vdots\\
\widetilde C_p
\end{pmatrix}.
\end{equation*}
where $\widetilde{C}$ denotes the matrix $C$ in controller coordinates. Even if the pair $(A,b)$ is controllable, the components of the output cannot in general be prescribed independently. The reason is that, in controller form, the state components satisfy differential relations. Consequently, different rows of $\widetilde C$ may impose differential constraints among the desired target components.

\begin{proposition}
Let $C_1,C_2\in\mathbb{R}^{1\times n}$, and suppose that, along every trajectory of the system, the corresponding outputs satisfy
\begin{equation*}
C_2x(t)=\frac{d}{dt}\big(C_1x(t)\big).
\end{equation*}
Then any pair of exactly output trackable targets $(f_1,f_2)$ must satisfy the compatibility condition $f_2(t)=f_1'(t)$ in the sense of distributions. If $f_1\in H^1(0,T)$, then necessarily $f_2=f_1'$ in $L^2(0,T)$.
\end{proposition}

\begin{proof}
If $C_1x=f_1$ and $C_2x=f_2$, then the assumed dynamical identity gives $f_2(t)=C_2x(t)=\frac{d}{dt}\big(C_1x(t)\big)=f_1'(t).$ This proves the compatibility condition.
\end{proof}

\begin{example}
Consider the system $x_1'(t)=x_2(t)$ and $x_2'(t)=-2x_1(t)-3x_2(t)+u(t)$, with full-state output $y_1(t)=x_1(t)$ and $y_2(t)=x_2(t)$. Then $y_1(t)=x_1(t),$ $y_2(t)=x_2(t)=x_1'(t).$ Therefore, simultaneous exact output tracking of $y_1(t)=f_1(t),$ $y_2(t)=f_2(t)$ is possible only if $f_2(t)=f_1'(t).$

Thus arbitrary pairs $(f_1,f_2)$ cannot be tracked exactly, even though the pair $(A,b)$ is controllable. This obstruction is not a failure of state controllability; it is a compatibility condition imposed by the chosen output map.
\end{example}

\begin{remark}
The preceding example illustrates the distinction between state controllability and output tracking. Controllability of $(A,b)$ allows one to steer the state at a final time, but exact output tracking imposes constraints on the whole time trajectory of $Cx(t)$. When several outputs are prescribed simultaneously, these constraints appear as differential relations among the target components.
\end{remark}

\subsection{Projection viewpoint when exact tracking fails}

When the full output $Cx$ cannot track an arbitrary target $f\in L^2(0,T;\mathbb{R}^p)$, one may still track suitable scalar
projections of the output. Let $\theta\in\mathbb{R}^p$ be a unit vector and define $z_\theta(t)=\theta^T Cx(t).$  The projected tracking problem consists in finding \(u\) such that $\theta^T Cx(t)=\theta^T f(t).$ This reduces the multi-output problem to a scalar-output problem with output matrix $C_\theta:=\theta^T C.$ If the corresponding relative degree is finite and the projected target $\theta^Tf$ has the required Sobolev regularity, then the scalar theory above applies.

From the dual HUM perspective, this projection corresponds to restricting the adjoint source to the direction $\theta$. The associated quadratic functional is
\begin{equation*}
J_\theta(g)=\frac12\int_0^T |b^T\phi(t)|^2\,dt - \langle g,\theta^T f\rangle,
\end{equation*}
where $\phi$ solves the adjoint equation with source $C^T\theta g$. Thus the HUM construction selects the component of the target that is observable through the pair $(b^T,\phi)$ in the chosen output direction.

\begin{example}
Consider the double-integrator system $x_1'(t)=x_2(t)$ and $x_2'(t)=u(t)$, with $x(0)=0$,  with full-state output $y_1(t)=x_1(t),$ $y_2(t)=x_2(t)$. Then every trajectory satisfies $y_1'(t)=y_2(t).$ Therefore, an arbitrary target $f=(f_1,f_2)$ cannot be tracked exactly unless $f_2=f_1'.$ However, for a fixed direction $\theta\in\mathbb{R}^2$, the projected
output $z_\theta(t)=\theta^T y(t)$ defines a scalar tracking problem. Whenever the corresponding projected
target $\theta^Tf$ satisfies the Sobolev regularity determined by the relative degree of $C_\theta=\theta^TC$, the projected component can be tracked by the scalar construction above.
\end{example}

This scalar-input analysis gives a concrete interpretation of the abstract results of Sections~\ref{sec:2}--\ref{sec:4}. Exact output tracking is governed not only by the controllability of the state equation, but also by the interaction between the control direction $b$, the dynamics $A$, and the chosen output map $C$. In controller coordinates, this interaction appears through relative degree and through differential compatibility conditions among output components.

\section{Numerical illustrations}\label{sec:6}

This section illustrates the main qualitative conclusions of the previous sections. The goal is not to implement exactly the HUM formula $u=\Lambda_C^*G_C^{-1}f$, but rather to solve a penalized tracking problem which provides a stable numerical approximation of the minimum-norm output tracking control.

We consider
\begin{equation}\label{eq:17}
\min_{u}
J_\alpha(u) = \frac12\int_0^T |u(t)|^2\,dt + \frac{\alpha}{2}\int_0^T |Cx(t)-f(t)|^2\,dt,
\end{equation}
subject to $x'(t)=Ax(t)+Bu(t)$, $x(0)=0.$ For large values of $\alpha$, the second term enforces an approximate tracking constraint, while the control cost regularizes the inversion. All experiments use $\alpha=10^4$. Thus, \eqref{eq:17} illustrates the penalized HUM viewpoint in situations where exact inversion is either regularity-limited or numerically ill-conditioned.

The numerical examples below are intended to illustrate three phenomena: the dependence of exact tracking on the relative degree of the chosen output; the loss of regularity induced by high relative degree; the deterioration of exact tracking requirements in semi-discrete PDE models. All computations are performed with CasADi and IPOPT, using a direct transcription scheme with $500$ time nodes, unless otherwise stated.

\subsection{Finite-dimensional ODE example}

We first consider the second-order system $x_1'(t)=x_2(t)$, $x_2'(t)=-2x_1(t)-3x_2(t)+u(t)$. The purpose of this example is to show how the choice of the output operator $C$ changes the relative degree and, consequently, the
regularity required of the reference trajectory.

For $C=(1\;0)$, the transfer function is $G(s)=\frac{1}{s^2+3s+2},$ so that the relative degree is $r=2$. Therefore, the explicit inversion formula requires $f\in H^2(0,T)$, and the corresponding tracking
control is formally given by $u(t)=f''(t)+3f'(t)+2f(t)$. This is consistent with Theorem~\ref{th:v1}: tracking the first state component requires two derivatives of the target.

\begin{figure}[H]
	\includegraphics[width=0.5\textwidth]{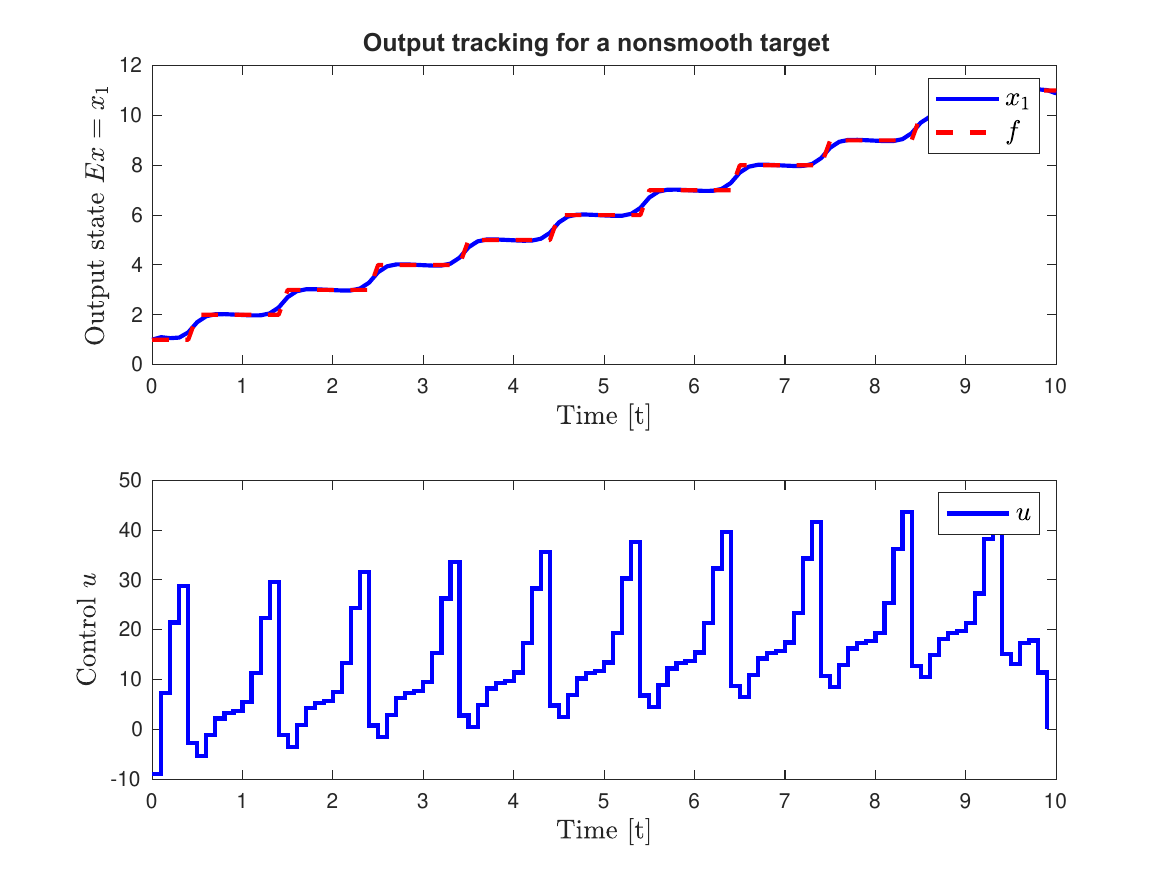}
\caption{Output $y(t)=x_1(t)$ (blue) and nonsmooth target $f(t)$ (red dashed). Despite insufficient regularity for classical tracking formulas, the variational approach provides a regularized/approximate tracking in $L^2$.}
	\label{fig4}
\end{figure}

For $C=(0\;1)$, the control appears after only one differentiation of
the output, and the relative degree is $r=1$. In this case, the target
only needs one derivative in $L^2(0,T)$, and the inversion formula
takes the form $u(t)=f'(t)-CAx(t)$, since $CB=1$. Numerically, this lower relative degree leads to smaller
tracking errors and smoother controls for smooth targets. We also consider the nonsmooth reference $f(t)=\lfloor t+1\rfloor$.  This function does not satisfy the $H^2$-regularity required by the
classical inversion formula for the output $C=(1\;0)$. Nevertheless,
the computation illustrates qualitatively that the penalized formulation can track nonsmooth references without explicit differentiation of the target (see Figure \ref{fig4}). 

\subsection{Semi-discrete wave equation}
We finally consider a semi-discrete model associated with the one-
dimensional wave equation with boundary control,
\begin{equation*}
\begin{cases}
w_{tt}-w_{xx}=0, & (x,t)\in (0,1)\times(0,T),\\
w(0,t)=u(t),\qquad w(1,t)=0, & t\in(0,T).
\end{cases}
\end{equation*}
The objective is to track the boundary observation $w_x(1,t)=f(t)$. After a centered finite-difference discretization with \(M\) interior points and mesh size \(h=(M+1)^{-1}\), the system becomes a finite-
dimensional linear control system of dimension \(2M\). The boundary tracking condition is approximated by $w_{M+1}(t)-w_M(t)=hf(t),$
 $w_{M+1}(t)=0,$ or equivalently $w_M(t)=-h f(t).$ Therefore, the tracking of the boundary derivative is reduced to the
tracking of the last discrete state component.

Applying the scalar-input analysis of Section~\ref{sec:5} shows that, for the semi-discrete system, exact tracking requires increasingly high Sobolev regularity of the target as the number of grid points increases. This explains the numerical ill-conditioning observed for fine discretizations: although the continuous wave equation admits boundary tracking in suitable regimes, the exact tracking requirement for the semi-discrete system becomes increasingly restrictive. This example illustrates a general warning: exact output tracking properties may deteriorate under spatial discretization, and robust numerical tracking should often be understood in an approximate or regularized sense.

\section{Further comments and open problems}\label{sec:8}
The results of this article show that tracking controllability is intrinsically a trajectory-space problem. Even for finite-dimensional systems, exact output tracking is governed by a trajectory-space observability inequality for the adjoint equation, involving derivative losses and structural restrictions on admissible targets. When this quantitative estimate fails, approximate output tracking remains meaningful and is characterized by a qualitative unique continuation property. 

In the scalar-input case, the controller canonical form makes these features explicit: compatibility conditions become differential relations among output components, while the Sobolev regularity required for exact tracking is determined by the relative degree. This viewpoint also suggests further directions beyond finite-dimensional systems, including tracking problems for PDE models and nonlinear settings \cite{DavronLissy2025, RisselTucsnak2023}. In particular, a sharper analysis of the output tracking observability inequality \eqref{eq:obs} for ODEs could clarify the dependence of constants and derivative losses, and may guide related sidewise or tracking estimates for PDEs \cite{asier2024}. Another relevant question is the simultaneous achievement of terminal state controllability and output tracking, namely enforcing $Cx(t)=f(t)$ on $[0,T]$ while steering the state to a prescribed final condition. This mixed objective involves additional compatibility constraints and deserves further study \cite{zuazua2022}.


\begin{thebibliography}{10}

\bibitem{AlbrechtGrasseWax1986}
F.~Albrecht, K.~A. Grasse, and N.~Wax.
\newblock Path controllability of linear input-output systems.
\newblock {\em IEEE Trans. Automat. Control}, 31:569--571, 1986.

\bibitem{asier2024}
J.~A. B\'arcena-Petisco and E.~Zuazua.
\newblock Tracking controllability for the heat equation.
\newblock {\em IEEE Trans. Automat. Control}, 70(3):1935--1940, 2025.

\bibitem{BrockettMesarovic1965}
R.~W. Brockett and M.~D. Mesarovi{\'c}.
\newblock The reproducibility of multivariable systems.
\newblock {\em Journal of Mathematical Analysis and Applications}, 11:548--563,
  1965.

\bibitem{darouach2025functional}
M.~Darouach and T.~Fernando.
\newblock On functional observability and functional observer design.
\newblock {\em Automatica J. IFAC}, 173:P. No. 112115, 9, 2025.

\bibitem{DavronLissy2025}
L.~Davron and P.~Lissy.
\newblock Exact output tracking for the one-dimensional heat equation and
  applications to the interpolation problem in {G}evrey classes of order 2,
  2025.

\bibitem{Fliess95}
M.~Fliess, J.~L\'evine, P.~Martin, and P.~Rouchon.
\newblock Flatness and defect of non-linear systems: introductory theory and
  examples.
\newblock {\em Internat. J. Control}, 61(6):1327--1361, 1995.

\bibitem{garcia2013alternative}
M.~I. Garc\'ia-Planas and J.~L. Dom\'inguez-Garc\'ia.
\newblock Alternative tests for functional and pointwise output-controllability
  of linear time-invariant systems.
\newblock {\em Systems Control Lett.}, 62(5):382--387, 2013.

\bibitem{GermaniMonaco83}
A.~Germani and S.~Monaco.
\newblock Functional output {$\varepsilon $}-controllability for linear systems
  on {H}ilbert spaces.
\newblock {\em Systems Control Lett.}, 2(5):313--320, 1983.

\bibitem{IsidoriBook}
A.~Isidori.
\newblock {\em Nonlinear control systems}.
\newblock Communications and Control Engineering Series. Springer-Verlag,
  Berlin, third edition, 1995.

\bibitem{kalman1960}
R.~E. Kalman.
\newblock Contributions to the theory of optimal control.
\newblock {\em Bol. Soc. Mat. Mexicana (2)}, 5:102--119, 1960.

\bibitem{lions1988controlabilite}
J.~L. Lions.
\newblock {\em Contr{\^o}labilit{\'e} exacte perturbations et stabilisation de
  syst{\`e}mes distribu{\'e}s. Tome 1, Contr{\^o}labilit{\'e} exacte.},
  volume~8.
\newblock Recherches en Mathematiques Appliqu{\'e}es, Masson, 1988.

\bibitem{RisselTucsnak2023}
M.~Rissel and M.~Tucsnak.
\newblock Approximate tracking controllability of systems with quadratic
  nonlinearities, 2025.

\bibitem{Ros70}
H.~H. Rosenbrock.
\newblock {\em State-space and multivariable theory}.
\newblock John Wiley \& Sons, Inc. [Wiley Interscience Division], New York,
  1970.

\bibitem{Rudinbook}
W.~Rudin.
\newblock {\em Functional analysis}.
\newblock McGraw-Hill, Inc., New York, second edition, 1991.

\bibitem{SainMassey69}
M.~K. Sain and J.~L. Massey.
\newblock Invertibility of linear time-invariant dynamical systems.
\newblock {\em IEEE Trans. Automatic Control}, AC-14:141--149, 1969.

\bibitem{zuazua2022}
Y.~Sara\c{c} and E.~Zuazua.
\newblock Sidewise profile control of 1-{D} waves.
\newblock {\em J. Optim. Theory Appl.}, 193(1-3):931--949, 2022.

\bibitem{Silverman1969}
L.~M. Silverman.
\newblock Inversion of multivariable linear systems.
\newblock {\em IEEE Trans. Automat. Control}, 14(3):270--276, 1969.

\bibitem{SilvermanPayne1971}
L.~M. Silverman and H.~J. Payne.
\newblock Input-output structure of linear systems with application to the
  decoupling problem.
\newblock {\em SIAM Journal on Control}, 9(2):199--233, 1971.

\bibitem{Wohltmann1985}
H.-W. Wohltmann.
\newblock Target path controllability of linear time-varying dynamical systems.
\newblock {\em IEEE Trans. Automat. Control}, 30:84--87, 1985.

\end{thebibliography}
\end{document}